\documentstyle{amsart}
\makeatletter
\def\LaTeX{\leavevmode L\raise.42ex
    \hbox{\kern-.3em\size{\sf@size}{0pt}\selectfont A}\kern-.15em\TeX}
\makeatother

\newcommand{\BibTeX}{{\rm B\kern-.05em{\sc i\kern-.025emb}\kern-.08em\TeX}}

\theoremstyle{plain}
\newtheorem{thm}{Theorem}[section]
\newtheorem{lem}[thm]{Lemma}
\newtheorem{cor}[thm]{Corollary}
\newtheorem{prop}[thm]{Proposition}

\theoremstyle{definition}
\newtheorem{defn}[thm]{Definition}

\theoremstyle{remark}
\newtheorem{rk}[thm]{Remark}
\newtheorem{notation}[thm]{Notation}
\newcommand{\thmref}[1]{Theorem~\ref{#1}}

\newcommand{\equref}[1]{{\rm ($\ref{#1}$)}}
\newcommand{\secref}[1]{\S\ref{#1}}
\newcommand{\lemref}[1]{Lemma~\ref{#1}}
\newcommand{\propref}[1]{Proposition~\ref{#1}}
\newcommand{\corref}[1]{Corollary~\ref{#1}}

\newcommand{\rkref}[1]{Remark~\ref{#1}}

\def\mapright#1{\smash{
   \mathop{\longrightarrow}\limits^{#1}}}
\def\mapleft#1{\smash{
   \mathop{\longleftarrow}\limits^{#1}}}
\def\mapdown#1{\Big\downarrow
     \rlap{$\vcenter{\hbox{$\scriptstyle#1$}}$}}
\def\exact#1#2#3{0\rightarrow{#1}\rightarrow{#2}\rightarrow{#3}\rightarrow0}
\begin{document}
\title{Variation of the Gieseker and Uhlenbeck Compactifications}
\author{Yi Hu$^1$}\footnote{Research partially support by NSF grant DMS
9401695}
\address{ Department of Mathematics, University of Michigan,  Ann Arbor, MI
48109}
\email{yihu@@math.lsa.umich.edu}
\thanks{}
\author{Wei-Ping Li}
\address{Department of Mathematics, Hong Kong University of Science and
Technology,
Clear Water Bay, Hong Kong}
\email{ mawpli@@uxmail.ust.hk}


\maketitle

\begin{abstract}
In this article, we study the variation of the Gieseker and
Uhlenbeck compactifications of the
moduli spaces of Mumford-Takemoto stable vector bundles of rank 2
by changing polarizations.
Some {\it canonical} rational morphisms among the Gieseker compactifications
are
proved to exist and their fibers are studied. As a consequence of studying the
morphisms
from the Gieseker compactifications to the Uhlebeck compactifications,
we show that there is an everywhere-defined {\it canonical} algebraic map
between two adjacent
Uhlenbeck compactifications which restricts to the identity on some
Zariski open subset.
\end{abstract}
\section{Introduction}\label{s:intro}

Let $X$ be an algebraic surface with $p_g=0$ and $H$ an ample divisor
over $X$. The moduli space ${\cal M}^\mu_H$ of
the Mumford-Takemoto $H$-stable rank two vector bundles
 has turned out to be a
key ingredient in the Donaldson theory of smooth topology of algebraic
surfaces. In
fact, Donaldson showed that the moduli space ${\cal N}_H$ of $SU(2)$-ASD
connections on $X$ with
respect to the Hodge metric induced by $H$ is homeomorphic to the moduli space
${\cal M}^\mu_H$
of Mumford-Takemoto $H$-stable rank two vector bundles. Hence the study of this
moduli
space is important for the application of the Donaldson theory.
It is obvious that the moduli space ${\cal M}^\mu_H \cong {\cal N}_H$ depends
on the polarization $H$.
The effect on the moduli space of $H$-stable bundles when changing the
polarization
has been considered before
by Donaldson \cite{D}, Friedman-Morgan \cite{FM}, Mong \cite{M} and Qin
\cite{Q1}, among others.
In particular, Qin \cite{Q1} gave a very systematic treatment.
\par

However, for many important applications, e.g., computing Donaldson's
polynomials,
just considering  the open variety ${\cal M}^\mu_H$ is not sufficient. In fact,
Donaldson
polynomials are computed on the Uhlenbeck compactification of ${\cal N}_H \cong
{\cal M}^\mu_H$.
So instead of considering variation of moduli spaces of Mumford-Takemoto
$H$-stable
rank-two vector bundles ${\cal M}^\mu_H$ for different $H$,
we take the step further to consider the variations of  the Gieseker and
Uhlenbeck
compactifications of the moduli space ${\cal M}^\mu_H$
of $H$-stable bundles.
\par
Gieseker constructed the moduli space ${\cal M}_H$ of $H$-semi-stable
torsion-free
coherent sheaves and showed that it is a projective scheme. Since ${\cal M}_H$
contains
the moduli space ${\cal M}^\mu_H$ as a Zariski open subset, ${\cal M}_H$ can be
considered
as a compactification of ${\cal M}^\mu_H$.  According to the Uhlenbeck weak
compactness theorem,
${\cal M}^\mu_H \cong {\cal N}_H$ also admits a Gauge theoretic
compactification.
This compactification is
called the Uhlenbeck compactification, denoted by $\overline{{\cal N}_H}$.

\par
It appears
that the works before  only considered variation
of ${\cal M}^\mu_H$ and set-theoretic comparison
of the moduli spaces ${\cal M}^\mu_H$ by varying $H$ (cf. \cite{Q1}).
However, in this paper, not only we take into account the variation of
compactifications of
${\cal M}^\mu_H$ but also our consideration is on the level of morphisms.
Namely, we  address
the existence of  morphisms amongst the moduli spaces. In particular,
we showed that there are {\it enough} canonical algebraic rational maps amongst
the
Gieseker compactifications and  canonical everywhere-defined algebraic maps
amongst
the Uhlenbeck compactifications. Moreover,  we gave some explicit description
of these morphisms and maps (see Theorem 5.1 and \S 7).
One of advantages of this  consideration is that
these maps carry considerable information which may allow one to trace
the geometry and topology from  one moduli space (and its compactifications)
to another.

\par
The main result of this paper may be summarized as follows.
Let $C_X$ be the K\"ahler cone of $X$ which is the closed convex cone in
$Num(X) \otimes {\Bbb R}$
spanned by all ample divisors.  There are certain natural wall and chamber
structures in
$C_X$ such that an ample divisor $H$ lies on a wall if and only if it possesses
non-universal
{\it strictly} MT $H$-semistable bundles. Let $C$ and $C'$ be two adjacent
chambers with a common face $F =
\overline{C} \cap \overline{C'}$. Pick up divisors $H$, $H'$, and $H_0$ in $C$,
$C'$
and $F$, respectively. Then there are two canonical rational morphisms
${\varphi}$ and
$\psi$ amongst the
Gieseker compactifications which descend to two everywhere-defined algebraic
maps
$\overline{\varphi}$ and $ \overline{\psi}$ amongst the Uhlenbeck
compactifications
$$\begin{matrix}{\cal M}_H(c_2)&\buildrel{\varphi}\over{-->}
& {\cal M}_{H_0}(c_2)& \buildrel{\psi}\over{<--} {\cal M}_{H'}(c_2)\\
\mapdown{\gamma_H}&&\mapdown{\gamma_{H_0}}&\mapdown{\gamma_{H'}}\\
\overline{{\cal N}_H}(c_2)&\mapright{\overline{\varphi}}&\overline{{\cal
N}}_{H_0}(c_2)
&\mapleft{\overline{\psi}} \overline{{\cal N}}_{H'}(c_2)
\end{matrix}$$
such that the above diagram commutes (see Theorem 5.1 and Theorem 7.8 for more
details).
 Here the morphisms $\gamma$ are the morphisms from the
Gieseker compactifications to their corresponding Uhlenbeck compactifications
as constructed
by J. Li \cite{Li}. Although $\varphi$ and $\phi$ are just rational maps
and hence are not surjective, $Im\varphi\cup Im\phi={\cal M}_{H_0}(c_2)$.
\par Another interesting result in this paper is about Uhlenbeck
compactifications. Uhlenbeck compactification $\overline {\cal N}_H(c_2)$ is,
in general, a closed subset of
$\coprod\limits_{j=0}^{c_2}\widetilde{\cal N}_H(j)\times Sym^{c_2-j}(X)$. It is
unknown whether $\overline {\cal N}_H(c_2)=\coprod\limits_{j=0}^{c_2}
\widetilde{\cal N}_H(j)\times Sym^{c_2-j}(X)$. When $p_g(X)=0$, we are able to
give an affirmative answer.
\par
Some of our considerations are inspired by a recent paper
of  Dolgachev and the first author
\cite{DH} where they treated the variational problem
of geometric invariant theory quotients.
However, we would like to point out that the variational problem of
the Gieseker compactifications and the Uhlenbeck compactifications
is considerably different from that of GIT!
Notably, the differences include, amongst other:\par
(1) in general, there are infinitely many moduli
spaces that are distinct to each other in  nature,
 while in the GIT case, the number of naturally distinct quotients is finite;
\par
(2) in general, there only exist rational maps among the Gieseker
compactifications, while
in the GIT case, morphisms among quotients are always defined everywhere.
Quite surprisingly, the maps among the Uhlenbeck compactifications are
defined everywhere.

\par
Another inspiration is Jun Li's paper on the relations between the Uhlenbeck
compactification
and the Gieseker compactification. Because a morphism from the Gieseker
compactification
to the Uhlenbeck compactification is constructed,
using results of the variation of the Gieseker compactification, we can
get some results on the variation of the Uhlenbeck compactification.
J. Morgan \cite{Mo}  also studied the map from Gieseker compactification
to Uhlenbeck compactification.

\par We mention that Friedman and Qin \cite{FQ} obtained stronger relations
among
the Gieseker compactifications and applied their results to good effect
on computing the Donaldson's invariants.  Also, after this work was completed,
 we learnt the work \cite{MW} and received
a copy of it.  However, neither of \cite{FQ} and \cite{MW}  stresses on
the Uhlenbeck compactifications.

\par
{\bf Acknowledgment}: We would like to thank the following people who
helped us in various ways during the course of this work: Robert Friedman,
Jun Li,  Zhenbo Qin, and Yungang Ye.
We also would like to thank Max-Planck-Institut f\"ur Mathematik
for inviting both of us to visit the institute in the summer of 1993.

\section{Background materials}
Let $X$ be an algebraic surface.
\begin{defn}\label{t:i1}
Let $V$ be a rank-two torsion-free coherent sheaf over $X$. Let $H$ be an ample
divisor on $X$
which will be called a stability polarization (or, a polarization for short).
$V$ is said to be Gieseker
$H$-stable ($H$-semi-stable) if for any rank one subsheaf $L$ of $V$,
$$\chi(L(nH))<\,\,(\le )\,\, {1\over 2} \chi(V(nH))\qquad \hbox{for $n\gg0$.}$$
$V$ is strictly Gieseker $H$-semi-stable if in addition there exists $L\subset
V$
such that
$$\chi(L(nH))= {1\over 2} \chi(V(nH)) \qquad \hbox{for $n\gg0$}.$$
\end{defn}
\par
There is another notion of stability
namely, the Mumford-Takemoto stability.
\begin{defn}\label{t:i2}
$V$ is said to be Mumford-Takemoto $H$-stable ($H$-semi-stable) if for any
rank one subsheaf $L$ of $V$,
$$c_1(L)\cdot H<\,\,(\le )\,\, {1\over 2} c_1(V)\cdot H.$$
$V$ is strictly Mumford-Takemoto $H$-semi-stable if in addition there exists
rank one subsheaf $L\subset V$ such that
$$c_1(L)\cdot H= {1\over 2} c_1(V)\cdot H.$$
\end{defn}

\begin{rk}\label{r:i1}
In this article, unless otherwise stated,
when we say $V$ is $H$-stable ($H$-semi-stable), we shall mean
Gieseker $H$-stable ($H$-semi-stable). We abbreviate
Mumford-Takemoto $H$-stable ($H$-semi-stable) as M-T$H$-stable
($H$-semi-stable).
\end{rk}
\par
Also, in this paper, the following convention will be adopted.
$V$, $V'$, etc.  represent  rank two torsion free coherent
sheaves and $L$, $L'$, $M$, $M'$,  ect.  represent  rank one torsion free
coherent
sheaves.

\par
Suppose $V$ is strictly $H$-semi-stable. Then
following Harder-Narishimhan filtration on semi-stable sheaves, we have that
$V$
sits in an exact sequence
$$\exact{L}{V}{L'}$$
with
$$\chi(L(nH))={1\over 2}\chi(V(nH)).$$
This exact sequence needs not to be unique but
$gr V=L\oplus L'$ is uniquely determined by $V$.
 We say that two strictly semi-stable bundles $V$ and $V'$ are $s$-equivalent
if $grV=grV'$
(see \cite{Gi}).

\par
Throughout this paper,  we use ${\cal M}_H(c_1, c_2)$, or ${\cal M}_H$ if
the  Chern classes are obvious from the context,
 to represent the moduli space of $H$-semi-stable sheaves
$V$ over $X$ with $c_1(V)=c_1$ and $c_2(V)=c_2$.
That is, ${\cal M}_H$ is the set of $H$-semi-stable sheaves modulo
$s$-equivalence.
Gieseker \cite{Gi} showed that ${\cal M}_H$
is a projective scheme.
We use ${\cal M}^\mu_H(c_1, c_2)$ (or ${\cal M}^\mu_H$) to represent M-T
$H$-stable vector bundles $V$ with
$c_1(V)=c_1$ and $c_2(V)=c_2$.

\section{Walls and Chambers}\label{s:chambers}
\begin{defn}
The K\"ahler cone $C_X$ of $X$ is the closed convex cone
in ${\rm Num}(X)\otimes {\Bbb R}$ spanned by ample divisors.
\end{defn}

For the purpose of comparing moduli spaces for varying polarizations,
we will introduce certain walls in the K\"ahler cone $C_X$. These walls
arise naturally from semi-stability.

Let $V$ be a rank 2 torsion-free coherent sheaf and $L$ be a subsheaf
of rank 1.
By Riemann-Roch formula, we have
$$\chi(V(nH))=\chi(V)+n^2H^2-nH\cdot K_X+nH\cdot c_1(V),$$
$$\chi(L(nH))=\chi(L)+{n^2\over 2}H-{n\over 2}H\cdot K_X+nH\cdot c_1(L).$$
\par
Hence
$$2\chi(L(nH))-\chi(V(nH))=(2\chi(L)-\chi(V))+n(2c_1(L)-c_1(V))\cdot H.$$
Therefore we obtain the following:

\par ($i$)
$V$ is $H$-stable if and only if for any given subsheaf $L$ one of the
following holds:
\par
{}~(1) $(2c_1(L)-c_1(V))\cdot H < 0$;
\par
{}~(2) $(2c_1(L)-c_1(V))\cdot H = 0$ but $2\chi(L)- \chi(V) <0$.
\par ($ii$)
Likewise, $V$ is strictly $H$-semi-stable if and only if for any given subsheaf
$L$ (1) or (2) of the
above holds except that  for some subsheaves $L$, we have
$(2c_1(L)-c_1(V))\cdot H = 0$ and $2\chi(L)- \chi(V) = 0$.

\par
In fact, in the above, we can always assume that the cokernel $V/L$ is torsion
free.
In particular, if $V$ is strictly $H$-semi-stable, $V$ sits in an exact
sequence
\begin{equation}\label{e:c1}
\exact{L}{V}{L'}
\end{equation}
with $(2c_1(L) -  c_1(V) )\cdot H = 0$ and $2\chi(L)=\chi(V)$.

Clearly,  that $V$ is  M-T $H$-stable implies that $V$ is $H$-stable. The
converse is not true, however.
Notice that $V$ is strictly M-T $H$-semi-stable if and only if for some
subsheaf $L$,
$(2c_1(L)-c_1(V))\cdot H = 0$, while in the Gieseker case we need to require
that
$2\chi(L)- \chi(V) = 0$. So the Gieseker stability is finer than M-T stability.
This is the main feature that distinguishes the variation problem
of Gieseker's stability from that of M-T stability.

\begin{defn}\label{t:c1}
Let $\tau\in {\rm Num}(X)$ be of the form
$2c_1(L) -  c_1$ where $L$ is a rank 1 sheaf. Assume further that
$-c\le \tau^2<0$ where $c$ is a fixed positive number.
We define the hyperplane of type $\tau$ as
$$W^{\tau}=\{\,h\in C_X |\,\tau\cdot h=0\}$$
$W$ is called a $c$-wall (or just a wall).
\end{defn}

\par
Let ${\cal W}$ be the set of $c$-walls in $C_X$.
It can be shown that for fixed $c$, the $c$-walls are locally finite.
Following  \cite{Q1}, we give the following definition.

\begin{defn}\label{t:c2}
A $c$-chamber (or just a chamber) $C$ in $C_X$ is a connected component of the
complement, $C_X - \bigcup_{W \in {\cal W}} W$, of the union of the $c$-walls.
A wall $W$ is called a wall of a chamber $C$ if $W \cap \overline{C}$ contains
a non-empty
open subset of $W$. In this case, The relative interior of
$W \cap \overline{C}$ is an  {\it open face} (or just {\it face}) of $C$. If
$F$ is a face of $C$,
then there is unique chamber $C' \ne C$ which also has $F$ as a face. In this
case,
the chamber $C$ and $C'$ lie on opposite sides of the wall containing the
common face $F$.
\end{defn}
It is obvious that each chamber (and each of its faces) is a convex cone in
$C_X$. In fact,
it is a polyhedral cone if its closure is contained entirely (except for the
origin)
in the interior of $C_X$.

\par
Now we fix $c=4c_2-c_1^2$, once and for all.
\par
Suppose $C$ and $C'$ share the same face $F$ lying on a wall  $W^{\tau}$. Then
$\tau\cdot H$ is either positive for all $H\in C$ or  negative
for all $H\in C$. A similar conclusion holds for $H'\in C'$.
Thus we may assume that $\tau \cdot C>0$ and $\tau\cdot C'<0$.
For simplicity, we shall
say that $C$ is the upper chamber and $C'$ is the lower chamber.
In many places of this paper,  we shall use
$H$ ($H'$) to represent an ample divisor in the chamber $C$ ($C'$)
respectively,
 $H_0$ to represent an ample divisor on the face $F$, and $\widetilde H$ to
represent
an arbitrary ample divisor.

\par
The following proposition partially justifies the definition of chambers.
\par
\begin{prop}\label{p:c1}
Let $C$ be a chamber in the K\"ahler cone. Then
${\cal M}_H={\cal M}_{H_1}$ for any two $H, H_1 \in C$.
\end{prop}
\begin{pf}
We shall prove this proposition by producing contradiction.
Without loss of generality, assume   that there exists $V\in {\cal
M}_H\backslash
{\cal M}_{H_1}$. Then there exists an exact sequence \equref{e:c1} such that
either we have $2c_1(L)\cdot H_1 > c_1(V)\cdot H_1$ or
we have $2c_1(L)\cdot H_1 = c_1(V)\cdot H_1$ but $2\chi(L)>\chi(V)$.
That is, if setting $\tau=2c_1(L)-c_1(V)$, then  $\tau\cdot H_1\ge 0$ and
\begin{equation}\label{e:c2}
\hbox{if }\tau\cdot H_1 = 0, \hbox{ then } 2\chi(L)>\chi(V).
\end{equation}
Since $V$ is $H$-semi-stable, we must have $\tau\cdot H\le 0$ and
\begin{equation}\label{e:c3}
\hbox{if }\tau\cdot H=0, \hbox{ then } 2\chi(L)\le \chi(V).
\end{equation}
Clearly  \equref{e:c2} and \equref{e:c3} cannot hold simultaneously.
Therefore $\tau$ cannot be  numerically trivial.
\par
Now we choose $H_2=(-\tau\cdot H)H_1+(\tau\cdot H_1)H$. Obviously, $\tau \cdot
H_2=0$.
Because the chamber $C$ is a convex cone,  $-\tau\cdot H\ge 0$ and
$ \tau\cdot H_1 \ge 0$,
 we obtain that $H_2$ is also an ample
line bundle in $C$. Since $\tau$ is not numerically
trivial, by Hodge index theorem, $\tau^2<0$.
\par
On the other hand, if we calculate the Chern classes from the exact sequence
\equref{e:c1}, we will get
$$ c_2(V)=c_1(L)\cdot(c_1(V)-c_1(L))+c_2(L)+c_2(L')
             \ge c_1(L)\cdot (c_1(V)-c_1(L)).$$
After some simplifications,  we get
$$\tau^2\ge c_1(V)^2-4c_2(V)=-c.$$
So $\tau$ defines a $c$-wall. Hence $H_2$ is in the chamber $C$ as well as on
the $c$-wall $W^\tau$, a contradiction.
\end{pf}

\section{Variation of ${\cal M}_H$ for different polarizations}\label{s:p}
Let $C$ and $C'$ be two chambers with a common face
$F \subset W^\tau$.
In this section, we will compare the moduli spaces ${\cal M}_H$, ${\cal
M}_{H'}$
and ${\cal M}_{H_0}$ where $H \in C, H' \in C'$, and $H_{0} \in F$.
\begin{defn}\label{t:p1}
$V$ is universally stable (semi-stable) if $V$ is stable (semi-stable) with
respect
to any polarization.
\end{defn}
In this section,  we will investigate what kind of  $H$-stable vector bundles
are not
$H_0$-semi-stable and so on. We will have a series of propositions of similar
nature.
\begin{prop}\label{p:p8}
Let $V$ be a $H$-semi-stable sheaf of rank 2.
Suppose that $V$ is not universally semi-stable.
Then $V$ must be $H$-stable. In another word, every semi-sable sheaves in
${\cal M}_H$ is $H$-stable
unless it is universally  semi-stable.
\end{prop}
\begin{pf}
Suppose $V$ is $H$-semi-stable, only two cases may happen: $V$ is M-T
$H$-stable,
or $V$ is  strictly M-T $H$-semi-stable.
If $V$ is M-T $H$-stable, then it must be  $H$-stable.
So assume that $V$ is strictly M-T $H$-semi-stable.
Then there exists a subsheaf
$L\subset V$ such that $V$ sits in the exact sequence \equref{e:c1} with
$2c_1(L)\cdot H=c_1(V)\cdot H$. If $\tau = 2c_1(L)-c_1(V)$ is not numerically
trivial,
then by Hodge index theorem,  $\tau^2<0$
and $\tau^2\ge c_1(V)^2-4c_2(V)=-c$. Hence $\tau$ defines a $c$-wall and $H$
lies
on the wall. This contradicts to the assumption that $H\in C$.
\par
Hence we must have $2c_1(L)-c_1(V)$ is numerically trivial. Therefore
$(2c_1(L)-c_1(V))
\cdot \tilde{H} =0$ for any ample divisor $\tilde{H}$, in particular, for the
ample divisor $H$.
 Since $V$ is $H$-semi-stable,
we must have $2\chi(L) \le \chi(V)$. From here, we  will show that $V$ is
$\tilde{H}$-
semi-stable for any $\tilde{H}$.
\par
In fact, assume that $M$ is a subsheaf of $V$. If $M$ is a subsheaf of $L$,
then
$c_1(M)\cdot \tilde{H}\le c_1(L)\cdot \tilde{H}$ and $\chi(M)\le \chi(L)\le
\displaystyle{1\over 2}\chi(V)$.
 Otherwise, $M$ admits an injection $M\hookrightarrow L'$. Either
$c_1(L')-c_1(M)$ is effective,
or $c_1(L')=c_1(M)$. In the first case, we obtain
$c_1(M)\cdot \tilde{H}<c_1(L)\cdot \tilde{H}
=\displaystyle{1\over 2}c_1(V)\cdot \widetilde H$.
In the latter case, i.e., $c_1(L')=c_1(M)$,
 we take double
dual of the exact sequence \equref{e:c1}, we get
\begin{equation}\label{e:p0}
0\rightarrow L^{**}\rightarrow V^{**}\rightarrow
                         L^{'**}I_Z\rightarrow 0,
\end{equation}
      $$V^{**}\hookleftarrow M^{**}=L^{'**}.$$
Hence $\ell(Z)=0$ and the exact sequence \equref{e:p0} splits, i.e.
$$V^{**}=L^{'**}\oplus L^{**}.$$
Therefore, the exact sequence \equref{e:c1} splits, i.e. $V=L \oplus L'$.
Since $2c_1(L)\cdot H=c_1(V)\cdot H=2c_1(L')\cdot H$ and $V$ is
$H$-semi-stable,
we must
have $2\chi(L)=\chi(V)=2\chi(L')$. Hence $2c_1(M)\cdot \tilde{H}=c_1(V)\cdot
\tilde{H}$ and
$2\chi(M)\le 2\chi(L')=\chi(V)$. This implies that  $V$ is
$\tilde{H}$-semi-stable. That is,
$V$ is a universally semi-stable sheaf.
But this contradicts to the assumption that
$V$ is not universally semi-stable.
\end{pf}

\begin{rk}\label{r:p2} The argument in  the proof above to show that the exact
sequence
\equref{e:c1} splits will be used (or referred) later on.
\end{rk}
\begin{cor}\label{p: c1}
Suppose $2c_1(L)-c_1(V)$ is numerically trivial, then any  non-splitting exact
sequence \equref{e:c1} gives a universally semi-stable sheaf provided
$2\chi(L)\le\chi(V)$.
\end{cor}

\begin{prop}\label{p:p2} Let  $V$ be strictly $\tilde{H}$-semi-stable sitting
in an
exact sequence \equref{e:c1} with $2c_1(L)\cdot \tilde{H} =c_1(V)\cdot
\tilde{H}$ and
$2\chi(L)=\chi(V)$.
\par ($i$)
If the exact sequence \equref{e:c1} doesn't split, then the subsheaf $L$
satisfying
$2c_1(L)\cdot \tilde{H} =c_1(V)\cdot \tilde{H}$ and $2\chi(L)=\chi(V)$ is
unique.
\par ($ii$)
Any $V$ sitting in \equref{e:c1} satisfying $2c_1(L)\cdot \tilde{H}
=c_1(V)\cdot \tilde{H}$
and $2\chi(L)=\chi(V)$  is strictly $\tilde{H}$-semi-stable.
\end{prop}
\begin{pf}
We only need to show the uniqueness of $L$.
\par
 Suppose otherwise, we have two exact sequences
$$\exact{L}{V}{L'},$$
$$\exact{M}{V}{M'}$$
satisfying $2c_1(L)\cdot \tilde{H}=2c_1(M)\cdot \tilde{H}=c_1(V)\cdot
\tilde{H}$,
$2\chi(L)=2\chi(M)=\chi(V)$ and $M\ne L$.
\par
If $M$ is a subsheaf of $L$, since $\chi(L)=\chi(M)$ and $c_1(L)\cdot
\tilde{H}=c_1(M)\cdot \tilde{H}$,
then $L=M$, a contradiction.
\par
Hence $M$ admits an injection into $L'$. By the similar argument as mentioned
in \rkref{r:p2}, the exact sequence splits: $V=L\oplus L'$ and $M=L'$. But we
have assume that
the sequence \equref{e:c1} does not split.
\end{pf}

{}From now on, we will be mainly concentrating on non-universally semi-stable
sheaves.

\begin{thm}\label{p:p3} Let $V$ be a non-universally semi-stable sheaf of rank
2.
Assume $V$ is $H_0$-stable.  Then one of the following holds.
\par  ($i$) If $V$ is M-T $H_0$-stable, then $V$ is M-T $H$-stable as well
as M-T $H'$-stable. In particular, $V$ is $H$-stable as well as  $H'$-stable.
\par  ($ii$) If $V$ is not M-T $H_0$-stable, then $V$ is either $H$-stable or
 $H'$-stable, but cannot be both.
\end{thm}
\begin{pf} Since $V$ is $H_0$-stable, then for subsheaf $L\subset V$ with
torsion
free cokernel, there exists
an exact sequence \equref{e:c1} such that either $2c_1(L)\cdot H_0< c_1(V)\cdot
H_0$
or $2c_1(L)\cdot H_0= c_1(V)\cdot H_0$ and $2\chi(L)<\chi(V)$.
\par
If $(2c_1(L)-c_1(V)\cdot H_0<0$, then $(2c_1(L)-c_1(V)\cdot H<0$ and
$(2c_1(L)-c_1(V)\cdot H'<0$. Otherwise there would exist a $c$-wall separating
$H_0$ with $H$ or $H_0$ with $H'$.
\par
If $2c_1(L)\cdot =c_1(V)\cdot H_0$, then $2\chi(L)<\chi(V)$. Since we assumed
that
$V$ is not universally stable, hence $2c_1(L)-c_1(V)$ is not numerically
trivial. Therefore $2c_1(L)-c_1(V)$ defines the  $c$-wall where $H_0$ lies.
Hence
$$\hbox{either}\qquad(2c_1(L)-c_1(V))\cdot H>0\qquad\hbox{ or }\qquad
(2c_1(L)-c_1(V))\cdot H<0.$$
\par
Since $V$ is $H_0$-stable, the exact sequence \equref{e:c1} doesn't split and
subsheaf $L$ satisfying $2c_1(L)\cdot =c_1(V)\cdot H_0$ is unique.
\par
Assume $(2c_1(L)-c_1(V))\cdot H<0$. For any subsheaf $M$ of $V$,
if $M$ is a subsheaf of $L$, we have
$2c_1(M)\cdot H<c_1(V)\cdot H$.
Otherwise, $M$ admits an injection into $L'$.
\par
 If $c_1(L')-c_1(M)$ is an
effective divisor, then $2c_1(M)\cdot H_0<2c_1(L')\cdot H_0=c_1(V)\cdot H_0$.
Hence
$(2c_1(M)-c_1(V))\cdot H_0<0$. Therefore $(2c_1(M)-c_1(V))\cdot H<0$.
\par If
$c_1(L')=c_1(M)$, then by the argument mentioned in \rkref{r:p2}, the exact
sequence \equref{e:c1} splits.  Since $V$ is $H_0$-stable, we get a
contradiction.
Hence $2c_1(M)\cdot H\cdot <c_1(V)\cdot H$. Therefore  $V$ is M-T $H$-stable.
\par Assume
$(2c_1(L)-c_1(V))\cdot H>0$, then $(2c_1(L)-c_1(V))\cdot H'<0$, hence by the
similar argument $V$ is
M-T $H'$-stable.
\par

The proof of ($i$) and ($ii$) will follow easily. For example, for ($ii$),
if $V$ is not  M-T $H_0$-stable, there will exist subsheaf $L\subset V$ such
that
$2c_1(L)\cdot H_0=c_1(V)\cdot H_0$. Hence if $2c_1(L)\cdot H<c_1(V)\cdot H$,
then
$2c_1(L)\cdot H'> c_1(V)\cdot H'$. In other words, if $V$ is $H$-stable, then
$V$
cannot be  $H'$ stable and vice-versa.
\end{pf}

\begin{thm}\label{p:p4} Let $V$ be a sheaf of rank 2 which is not universally
semi-stable.
Assume that $V$ is strictly $H_0$-semi-stable and sits in the non-splitting
exact sequence
\equref{e:c1}. Then $V$ is either $H$-stable or $H'$-stable, but can not be
both.
If the exact sequence \equref{e:c1} splits, then $V$ is neither $H$-stable nor
$H'$-stable.
\end{thm}
\begin{pf}
Assume that  $V$ is not $H$-stable nor $H'$-stable. Then there exist two exact
sequences
\begin{equation}\label{e:p1}
\exact{N}{V}{N'}
\end{equation}
\begin{equation}\label{e:p2}
\exact{M}{V}{M'}
\end{equation}
such that $2c_1(N)\cdot H\ge c_1(V)\cdot H$ and $2c_1(M)\cdot H'\ge c_1(V)\cdot
H'$.
Hence
\begin{equation}\label{e:p3}
2c_1(N)\cdot H_0\ge c_1(V)\cdot H_0\hbox{ and }2c_1(M)\cdot H_0\ge c_1(V)\cdot
H_0
\end{equation}
\par
Since $V$ is not universally semi-stable, it is easy to show that $2c_1(N)\cdot
H>c_1(V)\cdot H$ and $2c_1(M)\cdot H'>c_1(V)\cdot H'$.
Since $V$ is strictly $H_0$-semi-stable, we must have
\begin{equation}\label{e:p4}
2c_1(N)\cdot H_0\le c_1(V)\cdot H_0\hbox{ and }2c_1(M)\cdot H_0\le c_1(V)\cdot
H_0
\end{equation}
Combining \equref{e:p3} and \equref{e:p4}, we must have
$$2c_1(N)\cdot H_0=c_1(V)\cdot H_0=2c_1(M)\cdot H_0.$$
\par
Hence
\begin{equation}
2c_1(N)\cdot H'< c_1(V)\cdot H'\hbox{ and }2c_1(M)\cdot H< c_1(V)\cdot H
\end{equation}
Since $V$ is not universally stable, we have
\begin{equation}\label{e:p5}
2c_1(L)\cdot H\ne c_1(V)\cdot H\hbox{ and }2c_1(L)\cdot H'\ne c_1(V)\cdot H'
\end{equation}
If $2c_1(L)\cdot H<c_1(V)\cdot H$, then $N$ cannot be a subsheaf of $L$, since
otherwise, we would have $2c_1(L)\cdot H\ge 2c_1(N)\cdot H>c_1(V)\cdot H
2c_1(L)\cdot H$,
a contradiction. Hence $N$ admits an injection to $L'$. Then by the argument
mentioned in
\rkref{r:p2}, the exact sequence \equref{e:c1} splits, a contradiction.
\par
If $2c_1(L)\cdot H>c_1(V)\cdot H$, then $L$ cannot be a subsheaf of $M$, since
otherwise, we would have $c_1(V)\cdot H>2c_1(M)\cdot H\ge 2c_1(L)\cdot
H>c_1(V)\cdot H$,
a contradiction. Hence $L$ admits an injection to  $M'$. Then by the argument
mentioned in \rkref{r:p2}, the exact
sequence \equref{e:c1} splits, a contradiction.
\par Hence we proved that $V$ is either $H$-stable or $H'$-stable. It is easy
to
see that either $2c_1(L)\cdot H>c_1(V)\cdot H$ which implies that $V$ is not
$H$-stable, or $2c_1(L)\cdot H<c_1(V)\cdot H$, equivalently,
$2c_1(L)\cdot H'>c_1(V)\cdot H'$, which implies that $V$ is not $H'$-stable.
\par
If the exact sequence \equref{e:c1} splits, by the same argument as in the
previous paragraph, $V$ is neither $H$-stable nor $H'$-stable.
\end{pf}
Next, we give a criterion for strictly $H_0$-semi-stable sheaves.
\begin{prop}\label{p:p5} Assume that $V$ is not universally semi-stable.
$V$ is strictly $H_0$-semi-stable iff $V$ sits in the exact sequence
\equref{e:c1} where $c_1(L)\cdot H_0=c_1(L')\cdot H_0$ and $\chi(L)=\chi(L')$.
\end{prop}
\begin{pf} Easy.
\end{pf}

\begin{defn}\label{p:d2}
Let $L$, $L'$ be two rank one torsion free
coherent sheaves such that  $c_1(L)+c_1(L')=c_1$,
$c_1(L)\cdot c_1(L')+c_2(L)+c_2(L')=c_2$. If $c_1(L)-c_1(L')$ is numerically
trivial, then the pair ($L$, $L'$) is called $U$-pair. Otherwise, it is called
NU-pair.
\end{defn}

\begin{prop}\label{p:p6} Suppose ($L$, $L'$) is a NU-pair, $c_1(L)\cdot H_0
=c_1(L')\cdot H_0$ and $2c_1(L)\cdot H<c_1\cdot H$.
 Then
\par ($i$) every non-splitting exact sequence
\begin{equation}\label{e:p7}
\exact{L}{V}{L'}
\end{equation}
gives an $H$-stable sheaf $V$.
\par ($ii$) $V$ is $H_0$-stable if $\chi(L)<\chi(L')$;
\par ($iii$) $V$ is strictly $H_0$-semi-stable if $\chi(L)=\chi(L')$;
\par ($iv$) $V$ is $H_0$-unstable if $\chi(L)>\chi(L')$.
\end{prop}
\begin{pf} Consider a subsheaf $M\subset V$. If $M$ is a subsheaf of $L$, then
$$2c_1(M)\cdot H\le 2c_1(L)\cdot H< c_1(V)\cdot H.$$
Otherwise, $M$ admits an injection into $L'$. Hence
$$2c_1(M)\cdot H_0\le 2c_1(L')\cdot H_0=c_1(V)\cdot H_0.$$
\par
If $2c_1(M)\cdot H_0=c_1(V)\cdot H_0=2c_1(L')\cdot H_0$, by the argument
mentioned
in \rkref{r:p2}, the exact sequence \equref{e:p7} splits, a contradiction.
\par
Hence we must have $2c_1(M)\cdot H_0<c_1(V)\cdot H_0$. Hence $2c_1(M)\cdot
H<c_1(V)\cdot H$.
\par
The proves of (ii), (iii) and (iv) are quite straight forward.
\end{pf}

\begin{prop}\label{p:p7} Assume that $V$ is not universally semi-stable.
\par ($i$)
If $V$ is $H$-stable and strictly M-T $H_0$-semi-stable
sitting in the exact sequence \equref{e:c1} with $c_1(L)\cdot H_0=c_1(L')\cdot
H_0$,
then any sheaf $V'$ sitting in the non-splitting exact sequence
\begin{equation}\label{e:p6}
\exact{L'}{V'}{L}
\end{equation}
is $H'$-stable.
\par ($ii$)
If in addition, $V$ is not $H_0$-semi-stable, then $V'$ is $H_0$-semi-stable.
\end{prop}
\begin{pf} Since $V$ is not universally semi-stable, the pair ($L$, $L'$)
is an NU-pair. Because $V$ is $H$-stable and strictly M-T $H_0$-semi-stable,
 we obtain  $2c_1(L)\cdot H < c_1\cdot H$ and
$2c_1(L)\cdot H_0 = c_1 \cdot H_0$. Thus we get  $2c_1(L) \cdot H' > c_1\cdot
H'$, or equivalently,
$2c_1(L') \cdot H'< c_1 \cdot H'$. Apply \propref{p:p6} to the pair ($L'$,
$L$),
we get the conclusion.
\end{pf}

\begin{prop}\label{p:p9}
Fix the first Chern class $c_1$ and a $c$-wall $W^{\tau}$
where $c=4c_2-c_1^2$. For  $c_2'\ge c_2$, $W^{\tau}$ is also a $c'$-wall where
$c'=4c_2'-c_1^2$. Suppose that $\tau\cdot H_0=0$, $\tau\cdot H<0$ and
$\tau+c_1=2c_1(L)$ for some line bundle $L$. Then for $c_2'\gg0$, there
exists a $H$-stable sheaf $V$ with $c_1(V)=c_1$ and $c_2(V)=c_2'$ such that
it is not $H_0$-semi-stable.
\end{prop}
\begin{pf}
When $c_2(L')\gg 0$ , $Ext^1(L', L)\ne 0$. Hence there exists a  non-splitting
exact sequence
$$\exact{L}{V}{L'}$$
where $L'$ is a rank one subsheaf such that $c_1(L')=c_1-c_1(L)$ and
$c_2(L')=c_2'-c_2(L)-c_1(L)\cdot c_1(L')$. From \propref{p:p6}, $V$ is
$H$-stable with $c_2(V)=c_2'$ and $c_1(V)=c_1$.
\par
$$2\chi(L)=c_1^2(L)-c_1(L)\cdot K_X+2\chi({\cal O}_X),$$
$$\chi(V)=\displaystyle{c_1^2-c_1\cdot K_X\over 2}-c_2'+2\chi({\cal O}_X).$$
If $c_2'\gg 0$, then $2\chi(L)>\chi(V)$. Since $(2c_1(L)-c_1(V))\cdot H_0=0$,
hence $V$ is not  $H_0$-semi-stable.
\end{pf}

\par
In some propositions above, one may have noticed that we have used the term
``non-splitting  exact sequence'' several times. Suppose we have an exact
sequence
\equref{e:c1} with $c_1(L)\cdot H_0=c_1(L')\cdot H_0$ and $\chi(L)=\chi(L')$,
then
$V$ is strictly $H_0$-semi-stable and is $s$-equivalent to $gr(V)=L\oplus L'$.
 However, if $(L, L')$ is a NU-pair and if $V=L\oplus L'$, then $V$ is neither
$H$-stable nor $H'$-stable. Therefore if $Ext^1(L', L)=0$, then there would
exist a class in ${\cal M}_{H_0}$ represented by $V=L\oplus L'$ such that $V$
is
not s-equivalent to any $H_0$-semi-stable sheaf which  is  $H$-semi-stable,
nor is it s-equivalent to any $H_0$-semi-stable sheaf which is
$H'$-semi-stable.
In the following, we are going to show that
the situation above cannot happen.
This fact guarantees (ii) of  our main theorem in the next section.

\par
Let ${\cal F}$ and ${\cal F}'$ be two torsion free coherent sheaves.
Define (see \cite{Q2})
$$\chi({\cal F}',{\cal F})=\sum^2_{i=0}{\rm dim}Ext^i({\cal F}', {\cal F}).$$
\begin{prop}\label{p:n1}
$$\chi({\cal F}',{\cal F})=ch({\cal F}')^*\cdot ch({\cal F})\cdot td(X)_{H^4
      (X;{\Bbb Z})}$$
where $*$ acts on $H^{2i}(X;{\Bbb Z})$ by $(-1)^i\cdot Id$.
\end{prop}
\begin{cor}\label{c:n1}
Let $L$ and $L'$ be two torsion free rank one sheaves. Let $\tau=c_1(L)-c_1(L')
$. Let $V$ sit in the exact sequence \equref{e:c1}. Then
\begin{equation}\label{e: n1}
\chi(L', L)={\tau^2\over 4} -{K_X\cdot \tau\over 2} -c_2(V)+{c_1^2(V)\over 4}+
\chi({\cal O}_X)
\end{equation}
\end{cor}
\begin{pf}
{}From the exact sequence \equref{e:c1}, we get $\tau = 2c_1(L)-c_1(V)$ and
$c_2(L)+c_2(L')=c_2(V)-c_1(L)\cdot c_1(L')$.
$$\begin{array}{ll}
&\chi(L', L)=ch(L')^*\cdot ch(L)\cdot td(X)\\
=&(1-c_1(L')+\displaystyle {c_1(L')^2\over 2} -c_2(L'))\cdot (1+c_1(L)+
\displaystyle{ c_1(L)^2\over 2}
   -c_2(L))    \\
&~~\cdot (1-\displaystyle{ 1\over 2}K_X+\chi({\cal O}_X))_{H^4}\\
=&(1+c_1(L)-c_1(L')-c_1(L)\cdot c_1(L')+\displaystyle{ c_1(L)^2\over 2} +
\displaystyle{ c_1(L')^2\over 2}
      -c_2(L)-c_2(L'))\\
&~~\cdot (1-\displaystyle{ 1\over 2}K_X+\chi({\cal O}_X))_{H^4}\\
=&\displaystyle{ (c_1(L)-c_1(L'))^2\over 2}-
\displaystyle{ K_X\cdot (c_1(L)-c_1(L'))\over 2}-c_2(L)-c_2(L')
+\chi({\cal O}_X)\\
=&\displaystyle{ \tau^2\over 2}-\displaystyle{ K_X\cdot \tau\over 2}-c_2(V)+
(c_1(V)-c_1(L))\cdot c_1(L)+\chi({\cal O}_X)\\
=&\displaystyle{ \tau^2\over 2}-\displaystyle{ K_X\cdot \tau\over 2}
 -c_2(V)+\displaystyle{ c_1(V)^2-\tau^2\over 4}+\chi({\cal O}_X)\\
=&\displaystyle{ \tau^2\over 4}-\displaystyle{ K_X\cdot \tau\over 2}
 -\displaystyle{ 4c_2(V)-c_1(V)^2-3\chi({\cal O}_X)\over 4} +
\displaystyle{ \chi({\cal O}_X)\over 4}\\
=&\displaystyle{\tau^2\over 4} -\displaystyle{ K_X\cdot \tau\over 2} -
\displaystyle{ d\over 4}+\displaystyle{ \chi({\cal O}_X)\over 4}
\end{array}$$
where $d=4c_2(V)-c_1(V)^2-3\chi({\cal O}_X)$.
\end{pf}
\begin{prop}\label{p:n2}
With assumption on $L$ and $L'$ as in \propref{c:n1}. In addition, assume
that $\tau\cdot H_0=0$ and $d$, which is the virtual dimension of the moduli
space, is non-negative, then
$$-\chi(L, L')-\chi(L', L)>0.$$
\end{prop}
\begin{pf}Since ($L$, $L'$) is a NU-pair, $\tau$ is not numerically trivial,
 hence $\tau^2<0$.
{}From \corref{c:n1},
$$-\chi(L, L')-\chi(L', L)={d\over 2}+{-\tau^2-\chi({\cal O}_X)\over 2}
={d\over 2} +{-\tau^2-1+q\over 2}\ge {d\over 2}\ge 0.$$
Suppose $-\chi(L, L')-\chi(L', L)=0$, then $q=0, \tau^2=-1,
4c_2(V)-c_1^2(V)=3$.
\par
{}From the exact sequence \equref{e:c1}, we get
$$-\tau^2+4c_2(L)+4c_2(L')=4c_2(V)-c_1^2(V)=3.$$
Hence $1+4(c_2(L)+c_2(L'))=3$, or $4(c_2(L)+c_2(L'))=2$, impossible.
\end{pf}
\begin{cor}\label{c:n2}
${\rm dim}Ext^1(L, L')+{\rm dim}Ext^1(L', L)>0$. And there  either exists
a non-splitting exact sequence
$$\exact{L}{V}{L'},$$
or a  non-splitting exact sequence
$$\exact{L'}{V}{L}.$$
\end{cor}
\begin{pf}
Easy consequence of \propref{p:n2}.
\end{pf}
\begin{cor}\label{c:n3}
Suppose ($L$, $L'$) is a NU-pair satisfying  $c_1(L)\cdot H_0=c_1(L')\cdot
H_0$,
$\chi(L)=\chi(L')$,
then there exists a non-splitting
exact sequence
$$\exact{L}{V}{L'}\qquad\hbox{or}\qquad \exact{L'}{V}{L}$$
such that $V$ is strictly $H_0$-semi-stable.
\end{cor}
\begin{pf}
An easy consequence of  \corref{c:n2} and \propref{p:p5}.
\end{pf}

\section{Canonical rational morphisms among the Gieseker
compactifications}\label{s:g}

In this section, we shall draw some  conclusions on the variations of
Gieseker compactifications following many discussions in the previous section.

\par
Again, we place ourselves in the following situation.
Let $C$ and $C'$ be two chambers with a common face
$F \subset W^\tau$. We assume that
$C$ is the upper chamber and $C'$ is the lower chamber with respect to $\tau$.
That is, $C \cdot \tau > 0$ and $C' \cdot \tau <0$.
We will derive some canonical morphisms among the moduli spaces ${\cal M}_H$,
${\cal M}_{H'}$
and ${\cal M}_{H_0}$ where $H \in C, H' \in C'$, and $H_{0} \in F$.

\begin{thm}\label{t:n1}
There are two {\it canonical} rational algebraic maps
$${\cal M}_H  \buildrel{\varphi}\over{-->} {\cal M}_{H_0}
\buildrel{\psi}\over{<--} {\cal M}_{H'}$$
with the following properties:
\par ($i$) $\varphi$ and $\psi$ are isomorphisms over ${\cal M}^\mu_{H_0}
\subset {\cal M}_{H_0}$.
\par ($ii$) $Im\varphi \cup Im\psi={\cal M}_{H_0}$.
\par ($iii$) If $V\in {\cal M}_{H_0}$ is universally semi-stable, then the
inverse
            image $\varphi^{-1}(grV)$ or $\varphi^{-1}(V)$ consists of a single
point.
            The same conclusion holds for $\psi$.
\par($iv$)  If $V\in {\cal M}_{H_0}$ is not universally semi-stable and is
$H_0$-stable, then
             the inverse image $\varphi^{-1}(V)$ ($\psi^{-1}(V)$) consists of a
single point
             or is an empty set; More precisely,
\par~~~~($a$) If $V$ is MT $H_0$-stable, then $\varphi^{-1}(V)$
($\psi^{-1}(V)$) consists of
            a single point;
\par~~~~($b$)  If $V$ is strictly MT $H_0$-semistable, then $V$ sits in the
exact sequence  \equref{e:c1};
\par~~~~~~~~($b1$) If  $\chi(L)>\chi(L')$, $\varphi$ is not defined over ${\Bbb
P}(Ext^1(L, L')) \subset {\cal M}_H$,
           but $\psi$ sends ${\Bbb P}(Ext^1(L', L)) \subset {\cal M}_{H'}$
injectively to
           ${\cal M}_{H_0}$;
 \par~~~~~~~~($b2$) If $\chi(L) < \chi(L')$, $\psi$ is not defined over ${\Bbb
P}(Ext^1(L', L)) \subset {\cal M}_{H'}$,
           but $\varphi$ sends ${\Bbb P}(Ext^1(L, L')) \subset {\cal M}_H$
injectively to
           ${\cal M}_{H_0}$.
\par ($v$) If $V\in {\cal M}_{H_0}$ is not universally semi-stable and is
strictly $H_0$-semi-stable,
            then $V$ sits in the exact sequence \equref{e:c1} with
$\chi(L)=\chi(L')$,
             the inverse image of  $grV\in {\cal M}_{H_0}$ by $\varphi$ is
${\Bbb P}(Ext^1(L, L'))$,
            and the inverse image of  $grV\in {\cal M}_{H_0}$ by $\psi$ is
${\Bbb P}(Ext^1(L', L))$.
\end{thm}

\begin{pf}
First of all, for a given sheaf  $V$ in ${\cal M}_H$, if $V$ is also a
semi-stable
sheaf  with respect to $H_0$, then we can define a map which sends
$V \in {\cal M}_H$ (or $grV$) to $V$ (or $grV$) as a point in ${\cal M}_{H_0}$.
It is easy to see that this gives rise to a well-defined map
$\varphi$ from a Zariski open subset of ${\cal M}_H$ to ${\cal M}_{H_0}$.
Obviously, $\varphi$ is defined and restricts to the
identity over ${\cal M}^\mu_{H_0} \subset {\cal M}^\mu_H \subset {\cal M}_H $.
It remains to show the algebraicity of the map $\varphi$.
The proof is a standard one. So we only brief it.
Recall from the construction
of the moduli space ${\cal M}_H$ (see \cite{Gi}),
${\cal M}_H$ is the quotient of ${\cal Q}_H^{ss}$
by the group $PGL(N)$ (we adopt the notations from \cite{Li}).
 By the universality of the
quotient scheme, there is a universal quotient  sheaf ${\cal F}$
over $X \times {\cal Q}_H^{ss}$ with the usual
property. Now by the axiom of the coarse moduli,
there is a rational map from ${\cal Q}_H^{ss}$
to ${\cal M}_{H_0}$. Clearly this map respects the
group action (send an orbit of $PGL(N)$ to a point
by Proposition 4.2), thus by passing to the quotient,
we get a rational map from ${\cal M}_H$ to ${\cal M}_{H_0}$,
and this map is by definition the map $\varphi$.
Hence $\varphi$ is a morphism.

\par
The other map $\psi$ can be treated similarly.
\par
Property ($i$) and ($iii$) follows immediately from the above explanation.
\par
($iv $) and ($v$) follow as  consequences of \propref{p:p6} and \propref{p:p7}.

To prove ($ii$), take a $H_0$-semi-stable sheaf $V$.
\par
If $V$ is universally semi-stable, the conclusion follows by definition.
\par
If $V$ is not universally semi-stable but is  $H_0$-stable, then the conclusion
follows from \propref{p:p3}
\par
If $V$ is not universally semi-stable  and is strictly $H_0$-semi-stable, then
the conclusion follows from  \corref{c:n3} and \propref{p:p4}.
\end{pf}

\begin{prop}\label{v:p5} Let the situation be as in  \thmref{t:n1}.
If $c_2\gg0$, then the map $\varphi\colon {\cal M}_H(c_1, c_2)
 --> {\cal M}_{H_0}(c_1, c_2)$
and $\psi\colon {\cal M}_{H'}(c_1, c_2)  --> {\cal M}_{H_0}(c_1, c_2)$
are genuine rational maps (in the sense that they  can not be extended to
everywhere).
\end{prop}

\begin{pf}
It follows from  \propref{p:p9}.
\end{pf}

\section{From Gieseker's compactification to Uhlenbeck's
compactification}\label{s:u}
In this section, we will study the Uhlenbeck compactification of
 moduli spaces using its relation with the Gieseker compactification.
 We will use a technique established by Jun Li \cite{Li}
where he compared the Gieseker compactification and  the Uhlenbeck
compactification.
We assume that $q=0$ and $c_1=0$ through out this section.
Our analysis relies heavily on the results of Jun Li \cite{Li}.

\par
Following  the notations in \cite{Li},  let $H$ be an ample divisor and $g$ the
corresponding
Hodge metric on $X$. We use $\widetilde {\cal N}_{H}(j)$ to represent the
moduli
space of ASD connections, with respect to the Riemannian metric $g$,
 on an $SU(2)$ principal bundle $P$ over $X$ with
 $c_2(P)=j$,  and  ${\cal N}_{H}(j)$ to represent the moduli
space of irreducible ASD connections. ${\cal N}_{H}(j)$ is known by a
Donaldson's theorem
to be homeomorphic to the moduli
space  of Mumford-Takemoto $H$-stable vector bundles with $c_1=0$ and $c_2=j$.
We adopt the notation $\overline{\cal N}_H(c_2)$ to represent  the Uhlenbeck
compactification.
 The virtual dimension of the moduli space ${\cal M}_H(c_2)$ or $\overline{\cal
N}_H(c_2)$
is $d=4c_2-3$.

\par
Uhlenbeck compactification theorem tells us that $\overline{\cal N}_H(c_2)$ is
a closed subset of
$\coprod\limits_{j=0}^{c_2}\widetilde{\cal N}_H(j)\times Sym^{c_2-j}(X)$.
However, we didn't
know whether
$\overline {\cal N}_H(c_2)$ is the total space. The main conclusion
(\thmref{t:u2})
of this section will give a complete answer to this question.

\par
In what follows, we shall quote a useful theorem proved by J. Li (cf. Theorem
0.1, \cite{Li}).
\begin{thm}\label{t:u1} {\rm \cite{Li}}
There is a complex structure on $\overline{\cal N}_H(c_2)$ making it a reduced
projective scheme. Furthermore,
if we let ${\cal M}_H^{\mu}(c_2)$ be the open subset of ${\cal M}_H(c_2)$
consisting of locally free M-T $H$-stable sheaves and let
$\overline{\cal M}_H^{\mu}(c_2)$ be the closure of ${\cal M}_H^{\mu}(c_2)$
in ${\cal M}_H(c_2)$
endowed with reduced scheme structure, then there is a morphism
$$\gamma\colon \overline{\cal M}_H^{\mu}(c_2)\longrightarrow
\overline{\cal N}_H(c_2)$$
extending the homeomorphism between the set of M-T $H$-stable rank two vector
bundles and the set of gauge equivalent classes of irreducible ASD connections
with fixed Chern classes.
\end{thm}

It is known that when $c_2$ is large enough, ${\cal M}_H(c_2)$
is irreducible (\cite{GL}) and thus $\overline{\cal M}_H^{\mu}(c_2) = {\cal
M}_H(c_2)$
and $\gamma({\cal  M}_{H}(c_2)) = \gamma(\overline{\cal M}_H^{\mu}(c_2))$.
In \S 5 of \cite{Li}, a continuous map
$$\overline \sigma: \gamma(\overline{\cal M}_H^{\mu}(c_2)) \rightarrow
\coprod\limits^{c_2}_{j=0}\widetilde {\cal N}_{H}(j)\times Sym^{c_2-j}X$$ is
defined
and $\overline \sigma$ identifies  $\gamma(\overline{\cal M}_H^{\mu}(c_2))$
isomorphically with the Uhlenbeck compactification
$\overline{\cal N}_H(c_2) \subset \coprod\limits^{c_2}_{j=0}\widetilde {\cal
N}_{H}(j)\times
Sym^{c_2-j}X$.
We will use $\widetilde \sigma$ to stand for the map from $\gamma({\cal
M}_{H}(c_2))$
to $\coprod\limits^{c_2}_{j=0}\widetilde {\cal N}_{H}(j)\times Sym^{c_2-j}X$
defined in the
proof of Theorem 5 in \cite{Li}. If ${\cal  M}_{H}(c_2)$ is normal,
then this map $\widetilde \sigma$ is simply the map $\overline \sigma$.
Otherwise, $\widetilde \sigma$ is an extension of $\overline \sigma$.
\par
The definition of $\widetilde\sigma$ and $\overline\sigma$ are given in J. Li's
paper
\cite{Li}. It is recommended that the reader consult J. Li's paper to get
familiar with
these maps since in this and next section, we make use of these maps a lot.
\begin{rk}\label{r:u3}
In the rest of the paper, rather than directly working on  the Uhlenbeck
compactification
$\overline{\cal N}_H(c_2)$,
we will be  working on $\gamma({\cal  M}_{H}(c_2))$ instead.
One should keep in mind that $\gamma({\cal  M}_{H}(c_2))$ can be identified
via $\overline \sigma$ with the Uhlenbeck compactification $\overline{\cal
N}_H(c_2)$
when ${\cal  M}_{H}(c_2)$ is irreducible (and this can be ensured by requiring
$c_2$
to be large \cite{GL}). For small $c_2$, $\overline{\cal N}_H(c_2)$ is
contained
in $\gamma({\cal  M}_{H}(c_2))$ via the identification with
$\gamma(\overline{\cal M}_H^{\mu}(c_2))$.
In this case, $\gamma({\cal  M}_{H}(c_2))$ is slightly larger
than ${\cal  M}_{H}(c_2)$.
\end{rk}
\begin{notation}\label{n:u1}
Let $Z$ be a zero-cycle. $red(Z)$ will be the reduced scheme with
multiplicity counted at each point.
\end{notation}

We are going to prove the following proposition.
\begin{prop}\label{p:u1} Assume that $H$ is an ample divisor away from
$c$-walls.
\par{($i.$)} If $c_2\ge 2$, then ${\rm Im}\widetilde \sigma
=\coprod\limits^{c_2}_{j=0}{\cal N}_{H}(j)\times Sym^{c_2-j}X$.
\par {($ii.$)} If $c_2=1$, then ${\rm Im}\widetilde \sigma = {\cal N}_{H}(1)$.
In particular, $ {\cal N}_{H}(1)$ is compact.
\end{prop}
In order to prove the proposition, we divide  into several lemmas.
\begin{lem}\label{l:u1}
Suppose $V\in {\cal  M}_{H}(c_2)$ is strictly M-T $H$-semi-stable, then $V$
sits in the
exact sequence
\begin{equation}\label{e:u1}
\exact{I_Z}{V}{I_{Z'}}
\end{equation}
for some zero-cycles $Z$ and $Z'$ such that
$$\widetilde \sigma(\gamma(V))=({\cal O}_X\oplus {\cal O}_X,\, red(Z\cup
Z')).$$
\end{lem}
\begin{pf}
Since $V$ is strictly M-T $H$-semi-stable, hence there exist torsion free
coherent
sheaves of rank one $L$ and $L'$ such that $V$ sits in the exact sequence
$$\exact{L}{V}{L'}$$
with $c_1(L)\cdot H=0$.
\par Since $H$ is away from $c$-walls,  $c_1(L)$ is the trivial divisor by
Hodge index theorem. Hence we have the exact sequence \equref{e:u1}.
\par
By the proof of Lemma 3.3. in \cite{Li}, we get
$I_Z\oplus I_{Z'}\in \Gamma(\gamma(V))$. Therefore
$$\widetilde \sigma(\gamma(V))=({\cal O}_X\oplus {\cal O}_X,\, red(Z\cup Z'))
\in {\cal  M}_{H}(0)\times Sym^{c_2}(X).
$$
Note that ${\cal  M}_{H}(0)$ consists of a single point represented by
${\cal O}_X\oplus {\cal O}_X$.
\end{pf}
\begin{rk}\label{r:u2} Due to the same reason, the universally semi-stable
sheaves can only be sheaves sitting in the exact sequence \equref{e:u1}.
\end{rk}
\begin{lem}\label{l:u2}
Assume that $c_2\ge 2$. For any point
$({\cal O}_X\oplus {\cal O}_X, x)$ in ${\cal  M}_{H}(0)\times  Sym^{c_2}(X)$,
choose
a zero-cycle $Z$  of length $c_2$ such that $red(Z)=x$. Then there exists a
non-splitting exact sequence
\begin{equation}\label{e:u2}
\exact{I_Z}{V}{{\cal O}_X}
\end{equation}
such that $V$ is $H$-semi-stable.
\end{lem}
\begin{pf}
Let's calculate ${\rm dim} Ext^1({\cal O}_X, I_Z)=h^1(I_Z)$.
$$\begin{array}{ll}
h^1(I_Z)&=-\chi(I_Z)+h^0(I_Z)+h^2(I_Z)\\
& \ge -\chi(I_Z)=-(-c_2+1)\ge 1.\\
\end{array}$$
Hence there exists a nonsplitting exact sequence \equref{e:u2}.
\par
{\it Claim:} $V$ is $H$-semi-stable.
\par
In fact, let $M$ be a rank one subsheaf of $V$. Notice that
$$2\chi(I_Z)=2(-\ell(Z)+1)=-2\ell(Z)+2=-2c_2+2<-c_2+2=\chi(V).$$
If $M$ is a subsheaf  of $I_Z$, then $c_1(M)\cdot H\le 0$ and
$2\chi(M)\le 2\chi(I_Z)<\chi(V)$.
\par
Otherwise $M$ is a subsheaf of ${\cal O}_X$. Hence $c_1(M)\cdot H\le 0$. If
$c_1(M)\cdot H=0$, again $c_1(M)$ has to be the trivial divisor. Hence the
exact
consequence \equref{e:u2} splits, a contradiction. Therefore $c_1(M)\cdot H<0$.
\par
Thus,  we proved the claim.
\end{pf}
\begin{cor}\label{c:u1}
Assume that $H$ is an ample divisor away from  $c$-walls. Then
$${\rm Im}(\widetilde \sigma)\supset {\cal  M}_{H}(0)\times Sym^{c_2}(X).$$
\end{cor}
\begin{lem}\label{l:u3}
Given any element $(A, x)\in  \widetilde{\cal N}_{H}(j)\times Sym^{c_2-j}X$
where $j>0$, $A$
is an irreducible ASD corresponding to a M-T $H$-stable bundle $V_j$. In
particular
$${\cal N}_H(j)=\widetilde{\cal N}_H(j).$$
\par
Choose a zero-cycle $Z$ of length $\ell(Z)=c_2-j$ such that $red(Z)=x$.
Then the elementary transformation $V$ in
$$\exact{V}{V_j}{{\cal O}_Z}$$
is an M-T $H$-stable bundle in ${\cal  M}_{H}(c_2)$.
\end{lem}
\begin{pf}
Any reducible ASD in ${\cal N}_{H}(j)$ takes form of $L\oplus L^{-1}$ where $L$
is a line bundle
with $c_1(L)\cdot H=0$ and
$c_1(L)^2=-j<0$. Since $H$ is away from walls, such $L$ doesn't exist. Hence
$A$ corresponds to M-T $H$-stable bundle. The second statement is clear.
\end{pf}
\begin{cor}\label{c:u2}
With assumption and notations as in \lemref{l:u3}. Then
$\widetilde\sigma(V)=(A, x)$.
\end{cor}
\begin{pf}
$V^{**}=V_j$. By definition of $\widetilde\sigma$, $\widetilde\sigma(V)
=(A, x)$.
\end{pf}

Combination of \corref{c:u1} and \corref{c:u2} gives a proof of ($i$) of
\propref{p:u1}.
\begin{lem}\label{l:u4}
Assume $c_2=1$. Then there is no $H$-semi-stable sheaf $V$ which is strictly
M-T $H$-semi-stable.
\end{lem}
\begin{pf}
Suppose that $V$ is strictly M-T $H$-semi-stable, by \lemref{l:u1}, $V$ sits
in the exact sequence \equref{e:u1} with $\ell(Z)+\ell(Z')=c_2=1$.
Since $V$ is also $H$-semi-stable, $2\chi(I_Z)=-2\ell(Z)+2\le \chi(V)=1$.
Hence $\ell(Z)$ has to be one and $\ell(Z')=0$.
However,
$$\begin{array}{ll}
&{\rm dim}Ext^1({\cal O}_X, I_Z)\\
=&h^1(I_Z)=h^0(I_Z)+h^2(I_Z)-\chi(I_Z)\\
=&h^2({\cal O}_X)-(-\ell(Z)+1)=p_g=0.\\
\end{array}
$$
Hence the exact sequence (1) splits, i.e. $V={\cal O}_X\oplus I_Z$.
Since $2\chi({\cal O}_X)=2>\chi(V)=1$, $V$ will not be $H$-semi-stable, a
contradiction.
\end{pf}
This lemma proves ($ii$) of \propref{p:u1}.
Hence we finished the proof of  \propref{p:u1}.

\par
In the following, we consider the case where our polarization is on a face.
\begin{lem}\label{l:u5}
Suppose $H_0$ is an ample divisor on a face.
Suppose $V\in {\cal M}_{H_0}(c_2)$ is strictly M-T $H_0$-semi-stable, then
either $V$
sits in \equref{e:u1} and hence satisfies the conclusion of \lemref{l:u1},
or $V$ sits in the exact sequence
\begin{equation}\label{e:u3}
\exact{L\otimes I_Z}{V}{L^{-1}\otimes I_{Z'}}
\end{equation}
such that $c_1(L)\cdot H=0$ and $L$ is not the trivial line bundle. Moreover,
$$\widetilde\sigma(\gamma(V))=(L\oplus L^{-1},\,red(Z\cup Z')).$$
\end{lem}
\begin{pf}
The same as the proof  of \lemref{l:u1}
\end{pf}
For the first case in \lemref{l:u5}, \lemref{l:u1}, \lemref{l:u2} and
\lemref{l:u3}
still hold with some minor modifications. Let's prove the following lemma
which deals  with the latter case in \lemref{l:u5}.
\begin{lem}\label{l:u6} Without loss of generality, let's
assume $c_1(L)\cdot K_X\ge 0$. For any point $(L\oplus L^{-1}, x)$ in
$\widetilde{\cal N}_{H_0}^j\times Sym^{c_2-j}(X)$ where $0<j=-c_1(L)^2$, choose
a zero-cycle
$Z$ of length $c_2-j$ such that $red(Z)=x$. Then there exists an
$H_0$-semi-stable sheaf $V$ in the non-splitting exact sequence
\begin{equation}\label{e:u4}
\exact{L\otimes I_Z}{V}{L^{-1}}.
\end{equation}
\end{lem}
\begin{pf}
The proof is similar to that of \lemref{l:u2}.
$$\begin{array}{ll}
&{\rm dim}Ext^1(L^{-1}, L\otimes I_Z)=h^1( L^{\otimes 2}\otimes I_Z)\\
\ge&-\chi\bigl(L^{\otimes 2}\otimes I_Z\bigr)
 =-(\displaystyle{2c_1(L)\cdot (2c_1(L)-K_X)\over2}-\ell(Z)+1)\\
=&-2c_1^2(L)+c_1(L)\cdot K_X+\ell(Z)-1\\
=&c_2+c_1(L)\cdot K_X+j-1\ge c_2>0.\\
\end{array}$$
Hence there exists a non-splitting exact sequence \equref{e:u4}.

Notice that
$$\begin{array}{ll}
&~~2\chi(L\otimes I_Z)=2(\displaystyle{c_1(L)^2-c_1(L)\cdot K_X\over
2}-\ell(Z)+1)\\
&=-j-c_1(L)\cdot K_X-2\ell(Z)+2\\
&=-c_2+2-c_1(L)\cdot K_X-\ell(Z)\\
&\le -c_2+2=\chi(V).\\
\end{array}$$
\par
Let's prove that $V$ is $H_0$-semi-stable.
\par
Let $M$ be a rank one subsheaf of $V$.
If $M$ is a subsheaf of $L\otimes I_Z$, then
$$c_1(M)\cdot H_0\le c_1(L)\cdot H_0\qquad\hbox{and}\qquad
2\chi(M)\le 2\chi(L\otimes I_Z)\le \chi(V).$$
\par
Otherwise, $M$ is a subsheaf of $L^{-1}$. Hence $c_1(L^{-1})-c_1(M)$ is
effective
or trivial and $c_1(M)\cdot H_0\le c_1(L^{-1})\cdot H_0=0$. If $c_1(M)\cdot H_0
=c_1(L^{-1})\cdot H_0$, then $c_1(L^{-1})-c_1(M)$ is the  trivial divisor,
hence
the exact sequence \equref{e:u4} splits, a contradiction. Therefore,
$c_1(M)\cdot H_0<c_1(L^{-1})\cdot H_0=0$.
\par
Hence  $V$ is $H_0$-semi-stable.
\end{pf}
By the similar argument, we get the following proposition
\begin{prop}\label{p:u2} Assume that $H_0$ is an ample divisor on a face.
\par{($i.$)} If $c_2\ge 2$, then ${\rm Im}\widetilde \sigma
=\coprod\limits^{c_2}_{j=0}\widetilde {\cal N}_{H_0}(j)\times Sym^{c_2-j}X$.
\par {($ii.$)} If $c_2=1$, then ${\rm Im}\widetilde \sigma =\widetilde {\cal
N}_{H_0}(1)$.
In particular, $\widetilde {\cal N}_{H_0}(1)$ is compact.
\end{prop}
\begin{thm}\label{t:u2} Let $\widetilde H$ be an arbitrary polarization.
If ${\cal M}_{\widetilde H}(c_2)$ is irreducible, then
the Uhlenbeck
compactification is the total space, i.e.
$$\overline {\cal N}_{\widetilde H}(c_2)=\hbox{$\coprod\limits_{j=0}^{c_2}$}
\widetilde {\cal N}_{\widetilde H}(j)\times Sym^{c_2-j}(X)$$
when $c_2\ge 2$; When  $c_2 = 1$, we have
 $$\overline {\cal N}_{\widetilde H}(1)=
\widetilde {\cal N}_{\widetilde H}(1).$$
\end{thm}

\section{Canonical regular morphisms among the Uhlenbeck compactifications}
In the following, we are going to study the variation of the Uhlenbeck
compactifications.

\par
Let's recall some notations and results of J. Li \cite{Li}. Let ${\cal Q}_H$
be the Grothendieck's quotient scheme parameterizing all quotient sheaves
$F$ of ${\cal O}_X^{\oplus N}$ with $\hbox{det}F=H^{\oplus 2n}$ and
$c_2(F\otimes H^{-n})=c_2$. Let ${\cal Q}_H^{\mu}\subset {\cal Q}_H$ be the
open set consisting of all M-T $H$-semi-stable quotient sheaves. Let
${\cal Q}^{ss}_H\subset {\cal Q}_H$ be the open set consisting of all
$H$-semi-stable quotient sheaves.
\begin{rk}\label{r:u4}
J. Li constructed \cite{Li} a commutative diagram
$$\begin{matrix}{\cal Q}^{ss}_H\subset {\cal Q}_H^{\mu}
&\mapright{\gamma_{{\cal Q}_H}}
&{\Bbb P}^k\\
\mapdown{\pi}&&\Vert\\
{\cal M}_H(c_2)&\mapright{\gamma_H}&{\Bbb P}^k
\end{matrix}$$
for some projective space ${\Bbb P}^k$ such that (Lemma 3.2, \cite{Li})
$\gamma_{{\cal Q}_H}({\cal Q}^{\mu}_H\cap\overline{\cal Q}^{ss}_H)$ is
identical
to $\gamma_H({\cal M}_H(c_2))$ as sets where $\overline{\cal Q}^{ss}_H$ is the
closure of ${\cal Q}^{ss}_H$ in ${\cal Q}_H$.
\end{rk}

Since we are comparing spaces depending on different stability polarizations,
we will
use subscripts to distinguish different maps for different spaces, for example,
$\gamma_H$, $\overline\sigma_H$, ${\cal Q}_H$, etc.
\par
We also assume through out this section that the moduli spaces ${\cal
M}_H(c_2)$ and
 ${\cal M}_{H_0}(c_2)$ are normal. For example, when $c_2$ is sufficiently
 large,
 ${\cal M}_H(c_2)$ and ${\cal M}_{H_0}(c_2)$ are both normal, and generic
$H$-stable sheaves
are $H_0$-stable, too.
\par
Since we assumed that  ${\cal M}_H(c_2)$ and ${\cal M}_{H_0}(c_2)$ are both
normal, the map
$\widetilde\sigma$ defined and discussed in section \secref{s:u} becomes the
map $\overline
\sigma$ defined by J. Li.
\par
We are going to define a map $\overline\varphi$ from Uhlenbeck compactification
of
the moduli space $\overline{\cal N}_H(c_2)$ to $\overline{\cal N}_{H_0}(c_2)$.
\begin{notation}\label{n:h1}
We use $\ell_x(Q)$ to represent the length of the torsion sheaf $Q$ at
the point $x$.
\end{notation}
\par
Using Lemma 3.2 and Theorem 4 in [7], we can define a map
$$\overline\varphi\colon \overline {\cal N}_H(c_2)\longrightarrow
\overline {\cal N}_{H_0}(c_2)$$
as follows:
\par
Any element in $\overline {\cal N}_H(c_2)$ can be represented by
$\overline\sigma_H(\gamma_H(V))$ for some $H$-semi-stable sheaf $V\in {\cal
M}_H(c_2)$. We know that $V$ is either M-T $H_0$-stable or strictly
M-T $H_0$-semi-stable. If $V$ is the former, we define
$$\overline\varphi(\overline\sigma_H(\gamma_H(V)))=\overline\sigma_{H_0}(\gamma_{
H_0}(V)).$$
\par
If $V$ is the latter, then $V$ sits in an exact sequence
$$\exact{L\otimes I_Z}{V}{L^{-1}\otimes I_{Z'}}.$$
\par
Then we define
$$\overline\varphi(\overline\sigma_H(\gamma_H(V)))=(L\oplus L^{-1},
\hbox{$\sum$}(\ell_x(Z)x+\ell_x(Z')x)).$$
\begin{rk}\label{r:h2}
This map can be regarded as the induced map from $\varphi$ between Gieseker
compactifications.
\end{rk}
\begin{prop}\label{p:h1}
The map $\overline\varphi$ is well-defined.
\end{prop}

\noindent
{\it Proof}.
Suppose an element in $\overline {\cal N}_H(c_2)$ can be represented by
$\overline\sigma_H(\gamma_H(V))$ and $\overline\sigma_H(\gamma_H(V'))$ for
some $H$-semi-stable sheaves $V$ and $V'$. Then $V$ and $V'$ sit in the
exact sequences
$$\exact{V}{V^{**}}{Q},$$
$$\exact{V}{V^{'**}}{Q'}$$
where $Q$ and $Q'$ are supported at zero-dimensional schemes.
\par
By the definition of $\overline\sigma_H$ and $\gamma_H$ (see \cite{Li}),
$\overline\sigma_H(\gamma_H(V))
=(V^{**}, \hbox{$\sum$}\ell_x(Q)x)$
and
$\overline\sigma_H(\gamma_H(V'))=(V^{'**}, \hbox{$\sum$}\ell_x(Q')x)$.
Hence $V^{**}=V^{'**}$ and $\hbox{$\sum$}\ell_x(Q)x=\hbox{$\sum$}\ell_x(Q')x$.
\par
If $V$ is M-T $H_0$-stable, then $V^{**}=V^{'**}$ is M-T $H_0$-stable.
Hence $V'$ is $H_0$-stable. Then
$$\begin{array}{ll}
&\overline\varphi(\overline\sigma_H(\gamma_H(V))
=\overline\sigma_{H_0}(\gamma_{H_0}(V))\\
=&(V^{**}, \hbox{$\sum$}\ell_x(Q)x)
=(V^{'**}, \hbox{$\sum$}\ell_x(Q')x)\\
=&\overline\sigma_{H_0}(\gamma_{H_0}(V))
=\overline\varphi(\overline\sigma_H(\gamma_H(V')).\\
\end{array}$$
\par
Otherwise, $V^{**}=V^{'**}$ is strictly M-T $H_0$-semi-stable and $V$ sits
in the exact sequences
\begin{equation}\label{e:h0}
\exact{L\otimes I_{Z_1}}{V}{L^{-1}\otimes I_{Z_2}}
\end{equation}
and
$$\exact{V}{V^{**}}{Q}.$$
\par
By taking double dual of the exact sequence \equref{e:h0} , we get
$$\exact{L}{V^{**}=V^{'**}}{L^{-1}\otimes I_Z}.$$
\par
Hence we get
$$\exact{L/L\otimes I_{Z_1}}{V^{**}/V=Q}{L^{-1}\otimes I_Z/L^{-1}\otimes
I_{Z_2}}.$$
\par
Hence $$\ell_x(Q)+\ell_x(Z)=\ell_x(Z_1)+\ell_x(Z_2).$$
\par
Therefore
$$\begin{array}{ll}
&\overline\varphi(\overline\sigma_H(\gamma_H(V))\\
=&(L\oplus L^{-1}, \hbox{$\sum$}\ell_x(Z_1)x+\hbox{$\sum$}\ell_x(Z_2)x)\\
=&(L\oplus L^{-1}, \hbox{$\sum$}\ell_x(Q)x+\hbox{$\sum$}\ell_x(Z)x).\\
\end{array}$$
\par
By the similar argument as above for $V'$, we get
$$\begin{array}{ll}
&\overline\varphi(\overline\sigma_H(\gamma_H(V'))\\
=&(L\oplus L^{-1}, \hbox{$\sum$}\ell_x(Q')x+\hbox{$\sum$}\ell_x(Z)x).\\
\end{array}$$
Hence
$$\begin{array}{ll}
&\overline\varphi(\overline\sigma_H(\gamma_H(V))\\
=&(L\oplus L^{-1}, \hbox{$\sum$}\ell_x(Q)x+\hbox{$\sum$}\ell_x(Z)x)\\
=&(L\oplus L^{-1}, \hbox{$\sum$}\ell_x(Q')x+\hbox{$\sum$}\ell_x(Z)x)\\
=&\overline\varphi(\overline\sigma_H(\gamma_H(V')).\\
\end{array}$$  $\Box$

\begin{rk}\label{r:h3}
Since $\varphi$ is just a rational map, it might be expected that the induced
map
$\overline\varphi$ should also be only defined on a Zariski open subset.
However, the following two observations may be useful in understanding the
differences:
\par
(1) Uhlenbeck compactification losses track of Gieseker strictly
semi-stability.
It only respects M-T semi-stability ( see Lemma 3.3 in \cite{Li}).
\par
(2) When we regard the morphism defined by J. Li
$$\overline\sigma_H\circ\gamma_H\colon{\cal M}_H(c_2)\mapright{}\overline{\cal
N}_H(c_2)$$
as a blowing-down, then although
$${\cal M}_H(c_2)\buildrel{\varphi}\over{-->}{\cal M}_{H_0}(c_2)$$
is a rational map, after blowing downs on
${\cal M}_H(c_2)$ and ${\cal M}_{H_0}(c_2)$ respectively, the induced map
$\overline{\cal N}_H(c_2)\rightarrow \overline{\cal N}_{H_0}(c_2)$
becomes a well-defined map. Thus we have the following commutative diagram
$$\begin{matrix}{\cal M}_H(c_2)&\buildrel{\varphi}\over{-->}&{\cal
M}_{H_0}(c_2)\\
\mapdown{\gamma_H}&&\mapdown{\gamma_{H_0}}\\
\overline{{\cal N}_H}(c_2)&\mapright{\overline{\varphi}}&\overline{{\cal
N}}_{H_0}(c_2).
\end{matrix}$$
\end{rk}
\par
Next, we are going to show that the map $\overline\varphi$ is continuous in the
classical
complex topology.
\par
\begin{thm}\label{t:h1}
The map
$$\overline\varphi\colon \overline{\cal N}_H(c_2)\longrightarrow
\overline{\cal N}_{H_0}(c_2)$$
is continuous in analytic topology.
\end{thm}

\begin{pf}
 The argument pretty much follows the argument
in the proof of theorem 5 in \cite{Li}.
\par
Since $\overline{\cal N}_H$ and $\overline{\cal N}_{H_0}$ are both compact, it
suffices
to show that if $\lim s_n=s$ in $\overline{\cal N}_H$ and $\lim
\overline\varphi
(s_n)=t$ in $\overline{\cal N}_{H_0}$, then $\overline\varphi(s)=t$.
\par
Since $\gamma_H (\overline{\cal M}_H^\mu )= \overline{\cal N}_H$
is compact, ${\cal M}_H^{\mu}$ is dense in $ \overline{\cal M}_H^\mu$,
and generic $H$-stable sheaves are also $H_0$-stable, it suffices to show
the following statement: assume that $\{V_i\}$ is a sequence of $H$- and also
$H_0$-stable locally
free sheaves, $\lim V_i=V$ in ${\cal M}_H$,  and $\lim\overline\varphi
\circ \overline\sigma_H\circ \gamma_H(V_i)=t$ in $\overline{\cal N}_{H_0}$.
Then
$\overline\varphi(V)=t$.
\par
If $V$ is $H$-stable and $H_0$-semi-stable, clearly, the map $\overline\varphi
$ in the neighborhood of
 $V$ is induced from $\varphi$. Since $\varphi$ is continuous,
$\overline\varphi$ is continuous at $V$.
\par
Now suppose $V$ is $H$-stable and not $H_0$ semi-stable. It is clear that
$V$ is strictly MT $H_0$-semi-stable.
Since continuity  is a local problem, we can consider things locally.
In classical topology
there exists an open subset $U$ of ${\cal M}_H$ containing $V$ and a
universal sheaf ${\cal V}$ over $U\times X$ such that for any $u\in U$,
${\cal V}|_u$ represents $u$ in ${\cal M}_H$. Since every $H$-stable
sheaf is $H_0$-semi-stable, by Cor 1.4 in \cite{Gi}, we know that
there exists an integer $N$ such that $h^i({\cal V}|_u(NH_0))=0$ for $i\ge 1$
and
$H^0({\cal V}|_u(NH_0))$ generates ${\cal V}|_u(NH_0)$. By base change
theorem, we see that ${\cal V}(NH_0)|_{U'}$ is a quotient of ${\cal O}^{\oplus
r}_{
U'\times X}$ where $U'$ is an open subset of $U$ containing $V$ and
$r=h^0({\cal V}|_u(NH_0))$.
\par
By the universality of the quotient scheme ${\cal Q}_{H_0}$, we see that
there exists an analytic morphism
$$f\colon U'\longrightarrow {\cal Q}_{H_0}.$$
\par
Without loss of generality, we may assume that $V_i\in U'$.
Hence $$\lim_{i\to \infty}f(V_i)=f(V)\qquad \hbox{in ${\cal Q}_{H_0}$}.$$
\par
Therefore $f(V)\in {\cal Q}^{\mu}_{H_0}\cap\overline{\cal Q}^{ss}_{H_0}$.
\par
Since $\overline\varphi\circ\overline\sigma_H\circ\gamma_H(V_i)=
\overline\sigma_{H_0}\circ\gamma_{{\cal Q}_{H_0}}\circ f(V_i)$,
we have
$$\begin{array}{ll}
&\lim\overline\varphi\circ\overline\sigma_H\circ\gamma_H(V_i)=
    \lim\overline\sigma_{H_0}\circ\gamma_{{\cal Q}_{H_0}}\circ f(V_i)\\
=&\overline\sigma_{H_0}\circ\gamma_{{\cal Q}_{H_0}}\circ f(V)=t=
\overline\varphi\circ\overline\sigma_H\circ\gamma_H(V),\\
\end{array}$$
where the last equality comes from the definition of $\overline\varphi$.
\end{pf}

\par
Using the following lemma, the above theorem implies immediately the
algebraicity of
the map $\overline\varphi\colon \overline{\cal N}_H(c_2)\longrightarrow
\overline{\cal N}_{H_0}(c_2)$.

\begin{lem}
Let $X$ and $Y$ be two algebraic varieties. Let $U$ be a Zariski dense open
subset of $X$ with an algebraic morphism $\varphi_U\colon U\longrightarrow Y$
which
extends to $\varphi\colon X\longrightarrow Y$ continuously in analytic
topology. Then $\varphi$ is an
algebraic morphism.
\end{lem}

\begin{pf}
Consider the graph of the maps $\varphi$ and $\varphi_U$:
\begin{equation}
\rm{graph}(\varphi)\subset X\times Y,
\end{equation}
\begin{equation}
\rm{graph}(\varphi_U)\subset U\times X\subset X\times Y.
\end{equation}
Take the closure of $\rm{graph}(\varphi_U)$ inside $X\times Y$, call it
$\overline{\rm{graph}(\varphi_U)}$. Since $\varphi$ is continuous, it is easy
to see that $\overline{\rm{graph}(\varphi_U)}={\rm{graph}(\varphi)}$.
\par
In fact, for any element $(x, y)$ in $\overline{\rm{graph}(\varphi_U)}$,
there exists a sequence $(x_n, y_n) \in{\rm{graph}(\varphi_U)}$ such that
$(x_n, y_n)\rightarrow (x, y)$ as $n\rightarrow \infty$. Since
$\varphi(x_n)=y_n$,
\begin{equation}
y=\lim_{n\rightarrow \infty}y_n=\lim_{n\rightarrow\infty}\varphi(x_n)=
\varphi(\lim_{n\rightarrow\infty}x_n)=\varphi(x).
\end{equation}
Hence $\varphi\colon X\rightarrow Y$ can be regarded as a composition of
\begin{equation}
X\buildrel{\cong}\over\longrightarrow \rm{graph}\varphi =
\overline{\rm{graph}(\varphi_U)}
\hookrightarrow X\times Y
\buildrel{proj}\over\longrightarrow Y.
\end{equation}
Therefore, $\varphi$ must be an analytic morphism.
Since $X$ and $Y$ are all algebraic varieties,
$\varphi$ must be an algebraic morphism.
\end{pf}

\begin{cor} The map $\overline\varphi\colon \overline{\cal
N}_H(c_2)\longrightarrow
\overline{\cal N}_{H_0}(c_2)$ is algebraic.
\end{cor}

\par
In the rest of the section, we will study the inverse image of the map
$\overline\varphi$.
\begin{lem}\label{l:h1} Let $L$ and $L'$ be line bundles.
Suppose $V$ sits in a non-splitting exact sequence
$$\exact{L\otimes I_Z}{V}{L'\otimes I_{Z'}}$$
with $2c_1(L)\cdot H_0=c_1(V)\cdot H_0$. Then $V$ is strictly M-T
$H_0$-semi-stable.
\par
If in addition, $2c_1(L)\cdot H<c_1(V)\cdot H$, then $V$ is M-T $H$-stable.
\end{lem}
\begin{pf}
That $V$ is strictly M-T $H_0$-semi-stable is clear.
\par
Let $M$ be a rank one subsheaf of $V$.
If $M$ is a subsheaf of $L\otimes I_Z$, then
$$2c_1(M)\cdot H\le 2c_1(L)\cdot H<c_1(V)\cdot H.$$
Otherwise, $M$ is a subsheaf of $L'\otimes I_{Z'}$. Hence
$$2c_1(M)\cdot H_0\le 2c_1(L')\cdot H_0=c_1(V)\cdot H_0.$$
\par
If $2c_1(M)\cdot H_0<2c_1(L')\cdot H_0 = c_1(V)\cdot H_0$,
then $2c_1(M)\cdot H<c_1(V)\cdot H$, since otherwise, $2c_1(M)-c_1(V)$ would
define an $c$-wall between $H$ and $H_0$, a contradiction.
\par
If $2c_1(M)\cdot H_0= 2c_1(L')\cdot H_0$, then $L'=M$. Hence the exact sequence
splits, a contradiction.
\end{pf}
Now we divide the study of inverse image into several cases.

\par
(i) Suppose $(A, x)=({\cal O}_X\oplus {\cal O}_X, x)\in {\cal N}_{H_0}(0)\times
Sym^{c_2}(X)$ for $c_2\ge 2$. Then
$\overline\varphi^{-1}((A, x))=(A, x)\in {\cal N}_{H}(0)\times Sym^{c_2}(X)$.
Hence the inverse image of a point in the lowest stratum is just a single
point.
\par
(ii)
Suppose $(A, x)\in \coprod\limits_{j\ge 1}^{c_2}{\cal N}_{H_0}(j)\times
Sym^{c_2-j}(X)$, then $A$ corresponds to an M-T $H_0$-stable vector bundles
$V_j$. $V_j$ is also M-T $H$-stable. Then it is easy to see that
$$\overline\varphi^{-1}((A, x))=(A,x)\in\coprod_{j\ge 1}^{c_2}{\cal N}_H(j)
\times Sym^{c_2-j}(X).$$ Hence the inverse image of $\overline\varphi$
of a single point in $\coprod\limits_{j\ge 1}^{c_2}{\cal N}_H(j)
\times Sym^{c_2-j}(X)$ is also just a single point.
\par
(iii) Suppose $(A, x)=(L\oplus L^{-1}, x)$ with $L\cdot H_0=0$ and
$c_1(L)^2=-j$,
 i.e.
$$A\in \widetilde {\cal N}_{H_0}(j)-{\cal N}_{H_0}(j).$$
Without loss of generality, we assume that $L\cdot H<0$. Assume that $(A, x)=
\overline\varphi(\widetilde\sigma\gamma_H(V))$. By the way how the maps
$\widetilde \sigma$ and $\overline \varphi$ are defined, it is easy to see
that $V$
sits in the non-splitting exact sequence
\begin{equation}\label{e:h2}
\exact{L\otimes I_Z}{V}{L^{-1}\otimes I_{Z'}}.
\end{equation}
 By \lemref{l:h1}, $V$ is  M-T $H$-stable.

\par
The inverse image of $(A, x)$ is rather complicated due to the arbitrariness
 of $Z$ and $Z'$ in
the exact sequence \equref{e:h2}. We can only give a rough description of the
inverse image.
\par
Consider one extreme case where $V$ sits in a non-splitting exact sequence
$$\exact{L\otimes I_Z}{V}{L^{-1}}.$$
We know that $V$ is  $H$-stable. Clearly, $V^{**}$ sits in the exact sequence
$$\exact{L}{V^{**}}{L^{-1}}.$$

Hence
$$\overline\varphi^{-1}((L\oplus L^{-1}, x))\supset ({\Bbb P}(H^2(L^{\oplus
2}), x)
\subset {\cal M}_H(j)\times Sym^{c_2-j}(X).$$
\par
For other cases, it is easy to see that
$\overline\varphi^{-1}((L\oplus L^{-1}, x))$ consists of all
$$(V^{**}, x')\in {\cal M}_H(j')\times Sym^{c_2-j'}(X)$$
where $V$ sits in the exact sequence \equref{e:h2},  $c_2(V^{**})=j'$,
$x'\subset x$, and we
have the exact sequence
$\exact{V}{V^{**}}{{\cal O}_{Z''}}$ where $red(Z'')=x'$.
\par
In another word,
$\overline\varphi^{-1}((L\oplus L^{-1}, x))$ consists of
$$(V_{j'}, x')\in {\cal M}_H(j')\times Sym^{c_2-j'}(X)$$
where $c_2(V_{j'})=j'$, $x'\subset x$, $V_j$ is locally free sheave sitting in
the non-splitting exact sequence
$$\exact{L}{V_{j'}}{L^{-1}I_{Z'}}$$
 where $red(Z')=x-x'$.
\par
In general, the preimage of the map $\overline\varphi$ may contain points in
$\overline
{\cal N}_H(c_2)$ from every stratum ${\cal N}_H(j')$ for $j'\ge j$.
The intersection $\overline\varphi^{-1}((A, x))\cap{\cal N}_H(j')$
of the preimage with each stratum may not be closed. But put all these strata
together, $\overline\varphi^{-1}((A, x))$ will be closed.

\par
Summarize the above, we have
\begin{thm} Let the situation be as in Theorem 5.1. Then we have the following
commutative diagram
$$\begin{matrix}
{\cal M}_H(c_2)&\buildrel{\varphi}\over{-->}& {\cal M}_{H_0}(c_2)&
\buildrel{\psi}\over{<--} &{\cal M}_{H'}(c_2)\\
\mapdown{\gamma_H}&&\mapdown{\gamma_{H_0}}&&\mapdown{\gamma_{H'}}\\
\overline{{\cal N}_H}(c_2)&\mapright{\overline{\varphi}}&\overline{{\cal
N}}_{H_0}(c_2)
&\mapleft{\overline{\psi}} &\overline{{\cal N}}_{H'}(c_2)
\end{matrix}$$
such that
\par ($i$) $\overline{\varphi}$ and $\overline{\psi}$ are induced from
$\varphi$ and $\psi$;
both are well-defined everywhere and are  algebraic maps;
\par ($ii$) $\overline{\varphi}$ and $\overline{\psi}$ are homeomorphisms over
the Zariski open subset
$\coprod\limits_{j = 0}^{c_2}{\cal N}_{H_0}(j)
\times Sym^{c_2-j}(X)$;
\par ($iii$) Let $(A, x) \in (\widetilde {\cal N}_{H_0}(j)-{\cal N}_{H_0}(j))
\times Sym^{c_2-j}(X)$.
The the preimages of $\overline{\varphi}$ and $\overline{\psi}$ over $(A, x)$
are contained
in $\coprod\limits_{j' \ge j} \widetilde {\cal N}_{H_0}(j') \times
Sym^{c_2-j'}(X)$.
\end{thm}

\end{document}